\documentclass{amsart}

\usepackage{epsfig}

\newcommand{\cL}{{\mathcal L}}
\newcommand{\cT}{{\mathcal T}}
\newcommand{\bZ}{{\hspace{-0.3pt}\mathbb Z}\hspace{0.3pt}}

\newcommand{\lcr}{\raisebox{-5pt}{\mbox{}\hspace{1pt}
                  \epsfig{file=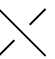}\hspace{1pt}\mbox{}}}
\newcommand{\ift}{\raisebox{-5pt}{\mbox{}\hspace{1pt}
                  \epsfig{file=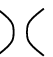}\hspace{1pt}\mbox{}}}
\newcommand{\zer}{\raisebox{-5pt}{\mbox{}\hspace{1pt}
                  \epsfig{file=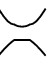}\hspace{1pt}\mbox{}}}
\newcommand{\trc}{\raisebox{-5pt}{\mbox{}\hspace{1pt}
                  \epsfig{file=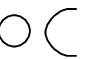}\hspace{1pt}\mbox{}}}
\newcommand{\lmk}{\raisebox{-5pt}{\mbox{}\hspace{1pt}
                  \epsfig{file=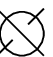}\hspace{1pt}\mbox{}}}
\newcommand{\mfr}{\raisebox{-5pt}{\mbox{}\hspace{1pt}
                  \epsfig{file=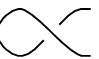}\hspace{1pt}\mbox{}}}
\newcommand{\pfr}{\raisebox{-5pt}{\mbox{}\hspace{1pt}
                  \epsfig{file=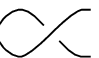}\hspace{1pt}\mbox{}}}
\newcommand{\nfr}{\raisebox{-5pt}{\mbox{}\hspace{1pt}
                  \epsfig{file=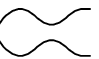}\hspace{1pt}\mbox{}}}
\newcommand{\mmk}{\raisebox{-5pt}{\mbox{}\hspace{1pt}
                  \epsfig{file=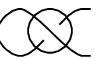}\hspace{1pt}\mbox{}}}
\newcommand{\pmk}{\raisebox{-5pt}{\mbox{}\hspace{1pt}
                  \epsfig{file=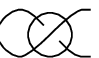}\hspace{1pt}\mbox{}}}
\newcommand{\cnl}{\raisebox{-5pt}{\mbox{}\hspace{1pt}
                  \epsfig{file=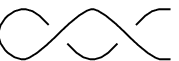}\hspace{1pt}\mbox{}}}
\newcommand{\cmm}{\raisebox{-5pt}{\mbox{}\hspace{1pt}
                  \epsfig{file=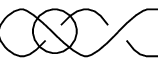}\hspace{1pt}\mbox{}}}
\newcommand{\cpm}{\raisebox{-5pt}{\mbox{}\hspace{1pt}
                  \epsfig{file=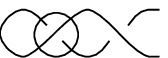}\hspace{1pt}\mbox{}}}
\newcommand{\tmk}{\raisebox{-12pt}{\mbox{}\hspace{1pt}
                  \epsfig{file=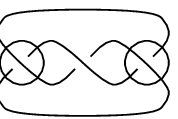}\hspace{1pt}\mbox{}}}
\newcommand{\trf}{\raisebox{-12pt}{\mbox{}\hspace{1pt}
                  \epsfig{file=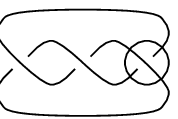}\hspace{1pt}\mbox{}}}
\newcommand{\frt}{\raisebox{18pt}{\mbox{}\hspace{1pt}
                  \epsfig{file=tref.eps,angle=180}\hspace{1pt}\mbox{}}}
\newcommand{\trz}{\raisebox{-12pt}{\mbox{}\hspace{1pt}
                  \epsfig{file=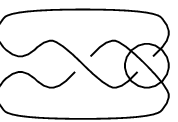}\hspace{1pt}\mbox{}}}
\newcommand{\zrt}{\raisebox{18pt}{\mbox{}\hspace{1pt}
                  \epsfig{file=tref-zero.eps,angle=180}\hspace{1pt}\mbox{}}}
\newcommand{\tri}{\raisebox{-12pt}{\mbox{}\hspace{1pt}
                  \epsfig{file=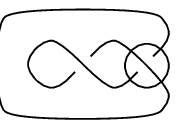}\hspace{1pt}\mbox{}}}
\newcommand{\irt}{\raisebox{18pt}{\mbox{}\hspace{1pt}
              \epsfig{file=tref-infinity.eps,angle=180}\hspace{1pt}\mbox{}}}

\newtheorem{lemma}{Lemma}

\newtheorem{prop}{Proposition}

\title{Skein Homology}

\author{Doug Bullock}

\address{Department of Mathematics, Boise State University, Boise, ID
  83725, USA\\ email:  bullock@@math.idbsu.edu}

\author{CHARLES FROHMAN}

\address{Department of Mathematics, University of Iowa, Iowa City, IA
  52245, USA\\ email: frohman@@math.uiowa.edu}

\author{JOANNA  KANIA-BARTOSZY\'NSKA}

\address{Department of Mathematics, Boise State University, Boise, ID
  83725, USA\\ email: kania@@math.idbsu.edu}

\begin{document}

\maketitle

\section{Introduction}

It is possible to introduce, for each skein module, a chain complex
whose homology is a $3$-manifold invariant and whose $0$-th level is
the original module.  In some sense the homology theories we
construct measure skein relations among skein relations, which mimics
Hilbert's theory of syzygies.

Skein modules were first introduced by Przytycki and Turaev
\cite{P1,T1}, with different motivations leading to various
constructions.  One in particular generalizes the Kauffman bracket
polynomial, itself a reformulation of the Jones polynomial.  We will
explicitly define the homology modules based on the Kauffman
bracket---a procedure easily generalized to other skein modules---and
demonstrate their nontriviality.

\section{The Kauffman Bracket Skein Module}

Generally, skein modules are constructed as follows. Take a free
module over some coefficient ring, spanned by isotopy classes of links
with some sort of decoration. Decorations include orientation of the
components and framings.  You then impose relations, which usually
come in two sorts: skein relations and framing relations.  For example,
to define the Kauffman bracket module we use the skein relation 
$\lcr= A\zer+A^{-1}\ift$, and the framing relation $L \cup\bigcirc
=-(A^2 + A^{-2})L$.  The first relation means that, outside of a small
embedded ball in $M$, the three links are identical, while inside they
appear as shown.  The second one means there is a trivial component of
the link that can be separated from the rest of the link by a ball.

Actually the Kauffman bracket skein module is built from framed links,
so you should think of the diagrams above as representing pieces of
flat annuli.  Let $\bZ[A,A^{-1}]$ be the ring of Laurent polynomials
with integer coefficients. Let $M$ be a $3$-manifold, and let
$\cL'(M)$ be the set of isotopy classes of framed links in $M$. A
framed link is an embedding of a disjoint union of annuli in $M$.  Let
$S'(M)$ be the smallest submodule of the free module
$\bZ[A,A^{-1}]\cL'(M)$ containing all sums of the form $\lcr-
A\zer-A^{-1}\ift$, and $L \cup\bigcirc +(A^2 + A^{-2})L$. The Kauffman
bracket skein module, $K(M)$, is the quotient
$\bZ[A,A^{-1}]\cL'(M)/S'(M)$.

One can define $K(M)$ beginning with a slightly different basis.  Let
$\cL(M)$ denote the set of framed links with no trivial components.
Let $S(M) < \bZ[A,A^{-1}]\cL(M)$ be the submodule of generated by all
sums of the form $\lcr- A\zer-A^{-1}\ift$ with the stipulation 
that you use $L \cup \bigcirc =-(A^2 + A^{-2})L$ to rewrite the
relation without trivial components.

\section{Kauffman Homology}

A crossing ball for a framed link is an embedding of the pair
$(B^3,D^2)$ so that, inside the crossing ball, the link is \lmk with
$D^2$ lying in the page.  A framed link with $n$ crossing balls is one
where the balls are ordered and disjoint. We number the crossing balls
from $1$ to $n$, corresponding to their ordering. Two such objects are
equivalent if there is an ambient isotopy of the links that carries
crossing ball to crossing ball in an order preserving fashion.  The
set of isotopy classes, excluding links with trivial components, will
be denoted by $\cL^n(M)$.

The $i$th ball operator,
\[ \partial_i : \bZ[A,A^{-1}]\cL^n(M) \rightarrow
\bZ[A,A^{-1}]\cL^{n-1}(M),\] 
is defined locally at the $i$th crossing ball by
\[\lmk \mapsto \lcr - A \zer - A^{-1} \ift,\]
along with any necessary applications of $L \cup \bigcirc =-(A^2 +
A^{-2})L$.  For example,
\begin{align*}
\partial_i\left(\pmk\right) &= \pfr - A\trc -A^{-1}\nfr \\
                 &= \pfr +A(A^2+A^{-2})\nfr - A^{-1}\nfr \\
                 &= \pfr +A^3\nfr.
\end{align*}
The ordering on the remaining crossing balls is induced by the
original one.  The boundary operator $\partial : \bZ[A,A^{-1}]\cL^n(M)
\rightarrow \bZ[A,A^{-1}]\cL^{n-1}(M)$ is the alternating sum $\sum_i
(-1)^i \partial_i$.  The boundary operator on framed links without
crossing balls is the zero map.  The proof that $\partial \circ
\partial =0$ is essentially the same as for singular homology.  Hence,
we have a chain complex, $(\bZ[A,A^{-1}]\cL^n(M), \partial)$.

The cycles are
\[Z_n(M)=
ker\{\partial: \bZ[A,A^{-1}]\cL^n(M) \rightarrow
\bZ[A,A^{-1}]\cL^{n-1}(M)\},\] the boundaries are
\[B_n(M)=
im\{\partial: \bZ[A,A^{-1}]\cL^{n+1}(M) \rightarrow
\bZ[A,A^{-1}]\cL^{n}(M)\},\] and the $n$th Kauffman homology of a
manifold is
\[K_n(M)= Z_n(M)/B_n(M).\]

We close this section with an example of a $2$-cycle in
$S^3$. Consider the element
\[ \cT = \tmk\]
of $\cL^2(S^3)$, with the crossing balls numbered from left to right. 
The ball operators are
\begin{align*}
\partial_1(\cT) &= \trf -A\trz - A^{-1}\tri, \quad \text{and}\\
\partial_2(\cT) &= \frt -A\zrt - A^{-1}\irt.
\end{align*}
 
A rotation of each diagram by $180^{\circ}$, shows that
$\partial_1(\cT)=\partial_2(\cT)$.  Hence,
$\partial(\cT)=-\partial_1(\cT)+\partial_2(\cT)=0$, so $\cT$ is a
$2$-cycle.

\section{Examples }

We will construct some cycles that represent nontrivial homology.  The
following technique for showing that a chain is not a boundary is due
to Chuck Livingston, improving on our original device.  Let $\zeta$ be
a primitive sixth root of unity.  Notice that this implies that $\zeta
+\zeta^{-1} = 1$. Since $\zeta^2$ is a third root of unity,
$\zeta^2+\zeta^{-2}=-1$. We define a map,
\[\epsilon : \bZ[A,A^{-1}]\cL^n \rightarrow \bZ[\zeta]\]
as follows. Each framed link $L$ has $\epsilon(L)=1$.  Set
$\epsilon(A)=\zeta$, $ \epsilon(A^{-1})=\zeta^{-1}$, and extend
linearly.

\begin{lemma}

$B_n \subseteq ker(\epsilon)$
\end{lemma}

\begin{proof}
Given $L\in \cL^{n+1}(M)$, 
\[\epsilon (\partial_i(L))= 
           \epsilon\left(\lcr - A \zer - A^{-1} \ift\right)=0,\]
provided there are no trivial components.  However, if a diagram
has a trivial component,
\[\epsilon (L \cup \bigcirc) =\epsilon(-(A^2 + A^{-2})L) = \epsilon(L),\]
so the computation above is still valid.

\end{proof}

The example from the last section  has $\epsilon(\cT) \neq 0$. Hence,
$\cT$ is not a boundary.
Therefore, $K_2(S^3) \neq 0$. The construction generalizes easily.

\begin{prop} Suppose that $L\in \cL^{2n}(M)$ and
that there exists a homeomorphism $h : M \rightarrow M$,
isotopic to the identity, which cyclically permutes the crossing balls.
Then $K_{2m}(M) \neq 0$. \end{prop}

We can apply this proposition to show that $K_{2n}(M) \neq 0$ for any
$3$-manifold $M$. For example, choose an annulus embedded in a ball in
$M$, with $2n$ kinks in it (all in the same direction) and surround
each kink with a crossing ball. There is an isotopy of $M$, fixed
outside a regular neighborhood of the knot, which permutes the kinks.
Another example is provided by embedding a $(2,2n)$ torus link into a
regular neighborhood of any knot in $M$.

How many different homology classes can be obtained via these examples?
That is, do examples built this way all turn out to be homologous, or
are there many distinct homology classes constructed the same way?

It turns out that framing forces first homology to be nontrivial.
\begin{prop}
For any $3$-manifold $M$, $K_1(M)\neq 0$.
\end{prop}

\begin{proof}
 Consider the one chain $\cpm - A^3\mmk.$ (The diagram can be
completed by any link in $M$.) Its boundary is
\[\cnl + A^3 \mfr -A^3\left(\mfr + A^{-3}\nfr\right) = 0,\]
so it's a cycle.  It's not a boundary for the usual reasons.
There's one with opposite framings, too: $\cmm - A^{-3}\pmk.$
\end{proof}

There are many open questions concerning skein homology. For instance,
what is $K_n(S^3)$? More generally, how is the topology of a $3$-manifold 
$M$ reflected in $K_n(M)$? 

We have described homology related to the Kauffman bracket skein
module. Obviously, the same construction could be carried out using
any other relations. The $0$-th level homology will recover the
original skein module.


\begin{thebibliography}{99}

\bibitem{P1} J. H. Przytycki, {\em Skein modules of 3-manifolds},
Bull.\ Pol.\ Acad.\ Sci.\ {\bf 39(1-2)} (1991) 91--100.
\bibitem{T1} V. G. Turaev, {\em The Conway and Kauffman modules of the
solid torus,} Zap.\ Nauchn.\ Sem.\ LOMI; English trans.\ in J. Soviet
Math. {\bf 167} (1988) 79-89.
\bibitem{kauffman1} L. Kauffman, {\em State models and the Jones polynomial},
Topology {\bf 26} no.\ 3 (1987) 395--401.
\end{thebibliography}
\end{document}